Article: CJB/2010/10

# A Construction to find any Circle through a Given Point

**Christopher J Bradley**

**Abstract:** A construction similar to Hagge's construction for circles through the orthocentre is shown to apply for any point.

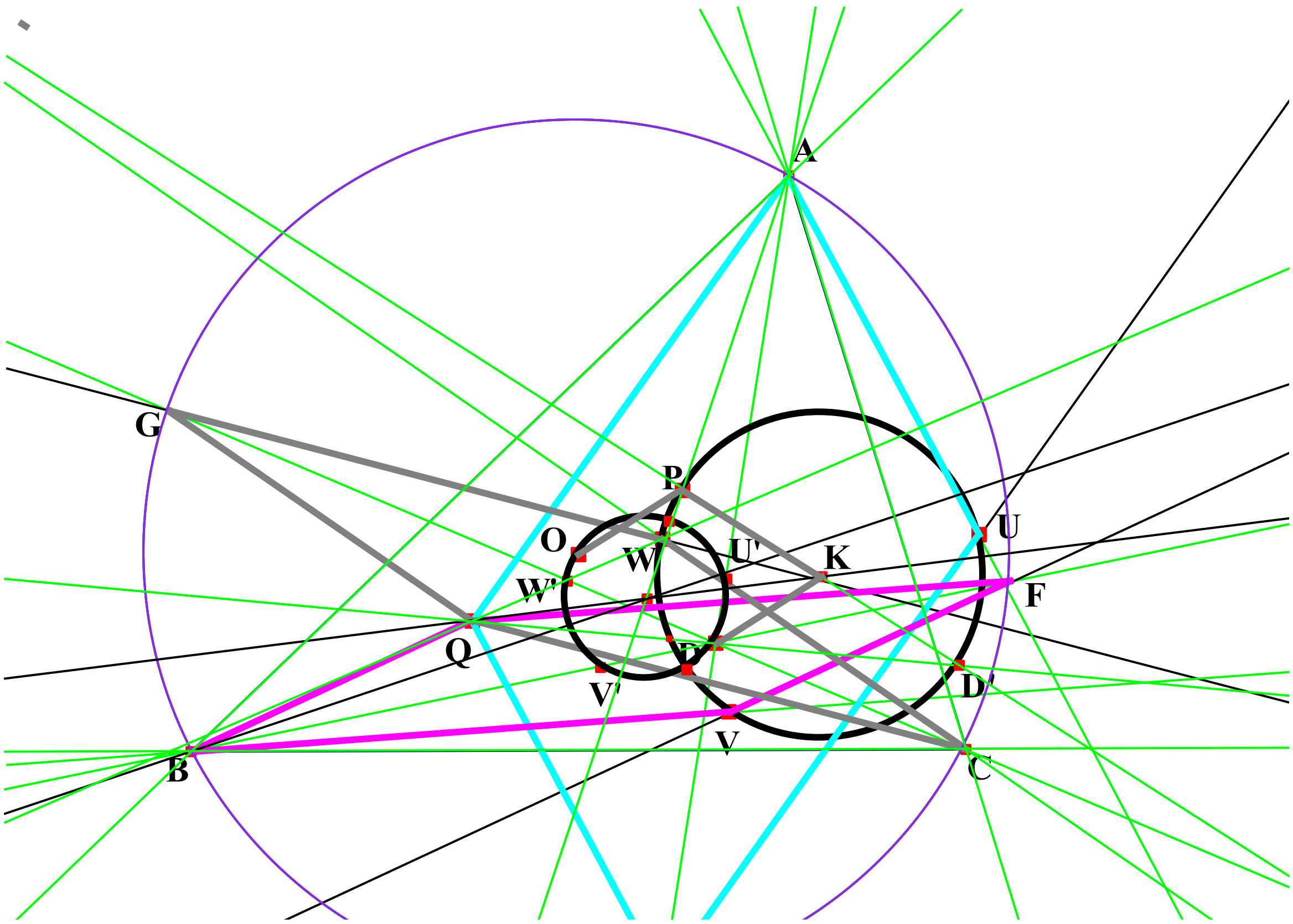

**Fig. 1**
**Showing the construction of a special circle through a point P**

## 1. Introduction

In 1907 Hagge wrote an article [1] in which he provided a construction that always provides a circle, now called a Hagge circle, through the orthocentre H of a triangle ABC. In the construction there is a key point D such that lines ADE, BDF, CDG meet the circumcircle at E, F, G respectively. And then E, F, G are reflected in the sides BC, CA, AB respectively to provide three points U, V, W that are concyclic with H. When you vary D you get a different circle through H. Later Peiser [2] showed that the isogonal conjugate of D, when rotated by 180° about the nine-point centre coincides with the centre of the Hagge circle. Furthermore triangles EFG and UVW are indirectly similar by a spiral similarity through the point D, as shown by the work of Speckman [3]. In the similarity triangle ABC is similar to a triangle XYZ inscribed in the Hagge circle, where triangle XYZ is in perspective to triangle UVW with vertex of perspective D. The disadvantage of the construction is the curious connection between the centre of the Hagge circle and the point D that generates it. The advantage is that because of the indirect similarity many interesting results, applications of the work by Speckman [3], exist for the configuration.

It has long been an objective that we should be able to provide a similar construction to provide *special circles* through any point P (not lying on the sides of ABC) but that has remained elusive until now. At the expense of the similarity we now have a method, which has the advantage that one can predetermine immediately the radius and centre of the resulting circle through P.

The method is as follows: If a special circle through P is required, then draw a segment POQ, where O is the circumcentre of triangle ABC and O is the midpoint of PQ. Then choose a point D so that OPKD is a parallelogram, where K is the desired centre of the circle and PK = OD is its desired radius. Then draw lines ADE, BDF, CDG where E, F, G lie on the circumcircle. Next construct the points U, V, W such that they are the fourth points of the parallelograms AQEU, BQFV, CQGW respectively. The key result is that circle UVW passes through P and has centre K and radius OD.

A second result is that if the diagonals of these three parallelograms meet at U', V', W' respectively, then circle U'V'W' is a circle passing through O and D, with OD as diameter, and with half the radius of the special circle through P. In fact QUU', QVV', QWW' are straight lines, so that the circles and the disposition of the points on them are directly similar by an enlargement factor two through Q (showing incidentally that one can construct similar circles through any point on the line QOP by choosing points on the diagonals QU, QV, QW suitably). All these result are illustrated in Fig. 1 above.

In the following sections we prove these results using Cartesian co-ordinates, taking the circumcircle to have radius 1 and centre at O. The working is technically quite elaborate, but the final results are beguilingly simple.

## 2. Position of the special circle and the points U, V, W

Take P to have co-ordinates (– k, 0), so that Q being the 180° rotation of P about O has co-ordinates (k, 0). Let D have co-ordinates (m, n). What we shall establish is that the radius of the special circle UVWP is $\sqrt{(m^2 + n^2)}$ and that its centre K has co-ordinates (m – k, n). Let A have co-ordinates $(2a/(1 + a^2), (1 - a^2)/(1 + a^2))$, and let B, C have corresponding parameters b, c. The equation of the line AD is

$$\{(n + 1)a^2 + (n - 1)\}x - \{(1 + a^2)m - 2a\}y + (1 - a^2)m - 2an = 0. \qquad (2.1)$$

This line meets the circumcircle with equation $x^2 + y^2 = 1$ at the point E with co-ordinates (x, y), where

$$x = 2\{m(n + 1)a^2 - (1 + m^2 - n^2)a + m(1 - n)\}/\{(m^2 + (n + 1)^2)a^2 - 4ma + m^2 + (1 - n)^2\},$$
$$y = \{(n + 1)^2 - m^2)a^2 - 4mna + m^2 - (1 - n)^2\}/\{(m^2 + (n + 1)^2)a^2 - 4ma + m^2 + (1 - n)^2\}. \qquad (2.2)$$

Having obtained E the co-ordinates of the displacement QE may be obtained. Bearing in mind that AU = QE, the co-ordinates of U may now be found and they are (x, y), where

$$x = (1/s)\{-(k(m^2 + (n + 1)^2) - 2m(n + 1))a^4 + 4(km + n(n + 1))a^3 - 2(k(m^2 + n^2 + 1) + 2m)a^2 + 4(km + n(n - 1))a - k(m^2 + n^2 - 2n + 1) - 2m(n - 1)\},$$

$$y = (1/s)\{-2(m^2a^4 + 2m(n - 1)a^3 - 4na^2 + 2m(n + 1)a - m^2)\}, \tag{2.3}$$

with

$$s = (1 + a^2)\{(m^2 + (n + 1)^2)a^2 - 4ma + m^2 + (n - 1)^2\}. \tag{2.4}$$

The co-ordinates of V and W may now be obtained by replacing the letter a in Equations (2.2) – (2.4) by letters b and c respectively.

## 3. The circle UVWP

The equation of the circle with centre (– g, – f) and radius $\sqrt{(g^2 + f^2 - t)}$ is

$$x^2 + y^2 + 2gx + 2fy + t = 0. \tag{3.1}$$

Inserting the co-ordinates of U, V, W we obtain three equations for f, g, t. When these are substituted back in Equation (2.5), we obtain the equation of circle UVW, which is

$$x^2 + y^2 + 2x(k - m) - 2mny + k(k - 2m) = 0. \tag{3.2}$$

It can now be seen that circle UVW passes through P(– k, 0), has centre K(m – k, n) and radius $\sqrt{(m^2 + n^2)}$, as required to be proved.

## 4. The points U'V'W' and the circle U'V'W'DO

Having obtained the co-ordinates of U, it is straightforward to obtain the co-ordinates of U', the midpoint of QU. Its co-ordinates (x, y) are independent of K, since U' is also the midpoint of AE. They are given by

$$x = (1/s)\{(n(1 + a^2) - (1 - a^2))(2an - m(1 - a^2))\},$$

$$y = (1/s)\{(2a - m(1 + a^2))(2an - m(1 - a^2))\}, \tag{4.1}$$

where s is given by Equation (2.4). The co-ordinates of V' and W' may be written down from Equation (4.1) by replacing the letter a by letters b and c respectively.

The equation of circle U'V'W' may now be obtained and is

$$x^2 + y^2 - mx - ny = 0. \tag{4.2}$$

This circle clearly passes through both O and D and has centre (½m, ½n). In fact OD is a diameter and its radius is $\frac{1}{2}\sqrt{(m^2 + n^2)}$.

## 5. Another circle through O

Finally by letting P, O, Q coincide at O and forming the parallelograms AOEU'', BOFV'', COGW'' the circle U''V''W'' passes through O and has centre D.

Flat 4,
Terrill Court,
12-14 Apsley Road,
BRISTOL BS8 2SP